\pdfoutput=1
\RequirePackage{ifpdf}
\ifpdf 
\documentclass[pdftex]{sigma}
\else
\documentclass{sigma}
\fi

\usepackage{bbm}
\usepackage{mathrsfs}
\numberwithin{equation}{section}
\newtheorem{Theorem}{Theorem}[section]
\newtheorem{Corollary}[Theorem]{Corollary}
\newtheorem{Lemma}[Theorem]{Lemma}
\newtheorem{Proposition}[Theorem]{Proposition}
\newtheorem{Conjecture}[Theorem]{Conjecture}
{ \theoremstyle{definition}
\newtheorem{Remark}[Theorem]{Remark}
}
\DeclareMathOperator{\Aut}{Aut}
\DeclareMathOperator{\End}{End}
\DeclareMathOperator{\Res}{Res}
\DeclareMathOperator{\Vir}{Vir}

\begin{document}

\allowdisplaybreaks

\renewcommand{\thefootnote}{$\star$}

\renewcommand{\PaperNumber}{019}

\FirstPageHeading

\ShortArticleName{Vertex Algebras $\mathcal{W}(p)^{A_m}$ and $\mathcal{W}(p)^{D_m}$ and Constant Term Identities}

\ArticleName{Vertex Algebras $\boldsymbol{\mathcal{W}(p)^{A_m}}$ and $\boldsymbol{\mathcal{W}(p)^{D_m}}$\\
and Constant Term Identities\footnote{This paper is a~contribution to the Special Issue on New Directions in Lie Theory.
The full collection is available at
\href{http://www.emis.de/journals/SIGMA/LieTheory2014.html}{http://www.emis.de/journals/SIGMA/LieTheory2014.html}}}

\AuthorNameForHeading{D.~Adamovi\'c, X.~Lin and A.~Milas}

\Author{Dra\v{z}en ADAMOVI\'C~$^\dag$, Xianzu LIN~$^\ddag$ and Antun MILAS~$^\S$}

\Address{$^\dag$~Department of Mathematics, University of Zagreb, Bijenicka 30, 10000 Zagreb, Croatia}
\EmailD{\href{mailto:adamovicl@math.hr}{adamovic@math.hr}}
\URLaddressD{\url{http://web.math.pmf.unizg.hr/~adamovic/}}

\Address{$^\ddag$~College of Mathematics and Computer Science, Fujian Normal University,\\
\hphantom{$^\ddag$}~Fuzhou, 350108, China} \EmailD{\href{mailto:linxianzu@126.com}{linxianzu@126.com}}

\Address{$^\S$~Department of Mathematics and Statistics, SUNY-Albany,\\
\hphantom{$^\S$}~1400 Washington Avenue, Albany 12222,USA}
\EmailD{\href{mailto:amilas@albany.edu}{amilas@albany.edu}}
\URLaddressD{\url{http://www.albany.edu/~am815139/}}

\ArticleDates{Received October 03, 2014, in f\/inal form February 25, 2015; Published online March 05, 2015}

\Abstract{We consider $AD$-type orbifolds of the triplet vertex algebras $\mathcal{W}(p)$ extending the well-known $c=1$
orbifolds of lattice vertex algebras.
We study the structure of Zhu's algebras $A(\mathcal{W}(p)^{A_m})$ and $A(\mathcal{W}(p)^{D_m})$, where $A_m$ and $D_m$
are cyclic and dihedral groups, respectively.
A~combinatorial algorithm for classif\/ication of irreducible $\mathcal{W}(p)^\Gamma$-modules is developed, which relies
on a~family of constant term identities and properties of certain polynomials based on constant terms.
All these properties can be checked for small values of~$m$ and~$p$ with a~computer software.
As a~result, we argue that if certain constant term properties hold, the irreducible modules constructed
in~[\textit{Commun. Contemp. Math.} \textbf{15} (2013), 1350028, 30~pages; \textit{Internat.~J. Math.}
  \textbf{25} (2014), 1450001, 34~pages]
provide a~complete list of irreducible $\mathcal{W}(p)^{A_m}$ and $\mathcal{W}(p)^{D_m}$-modules.
This paper is a~continuation of our previous work on the $ADE$ subalgebras of the triplet vertex algebra
$\mathcal{W}(p)$.}


\Classification{17B69}

\renewcommand{\thefootnote}{\arabic{footnote}}
\setcounter{footnote}{0}

\vspace{-1mm}

\section{Introduction and notation}

Orbifold theory of vertex algebras has an interesting and rich development.
Although orbifolds have been studied even earlier in string theory, their f\/irst mathematical treatment goes back to
already classical work by Frenkel, Lepowsky and Meurman on the Moonshine module, which is constructed as
a~$\mathbb{Z}_2$-orbifold~\cite{flm}.
We should also mention, closely related to our investigation, an important construction of $c=1$ orbifolds~by
Ginsparg~\cite{gin}.
In the language of vertex algebra, Ginsparg was basically considering what we now call $ADE$-type orbifolds of the rank
one lattice vertex algebra of central charge one and associated orbifold characters and the partition function.
This line of work was later brought to even f\/irmer footing in the VOA literature by Dong, Griess and others
(see~\cite{dg} and references therein).
Our current line of work is an attempt to lift these classical results from rational to the setup of irrational vertex
algebras with central charge $1-\frac{6(p-1)}{p^2}$, $p \geq 2$.

Arguably, the most famous constant term identity is the one due to F.~Dyson who discovered it in connection to what we
now call the circular ensembles model in random matrix theory.
His conjectural identity, later proved by Gunson, Wilson and others, concerned the constant term of
\begin{gather*}
\prod\limits_{1 \leq i \neq j \leq n} \left(1-\frac{x_i}{x_j}\right)^{p},
\qquad
p \in \mathbb{N},
\end{gather*}
which can be elegantly expressed as a~ratio of factorials; observe that the rational function is up to a~monomial term
just the~$p$-th power of the discriminant.
Dyson's identities are later generalized by Morris, Kadell, Aomoto and others by adding extra terms to the discriminant
part.
For example, one includes additional variable $x_0$ (e.g.
Morris and Aomoto's identities) or perhaps symmetric functions in $x_i$ (as in Kadell's identitiy).
All these identities are in a~way related to Selberg's integral (in fact, Morris' identity seems to be {equivalent} to
it).
For a~good review of this subject see~\cite{fw}.

Constant term identities are also related to computation of correlation functions arising from vertex operators in
$2$-dimensional CFT.
This can be either done in the integral form leading to Selberg integral, or purely formally by using constant terms of
generating functions~-- the method that we advertise here.
There are three main computational ingredients that lead to aforementioned expressions and eventually to considerations
of constant term identities
\begin{itemize}\itemsep=0pt
\item Various calculations with products of bosonic vertex operators $Y(e^{\alpha_i},x_1) \cdots Y(e^{\alpha_n},x_n)$ and
their normal ordering.
Already such considerations lead to Dyson-type expressions (see for instance~\cite{am-1}).

\item Action of certain distinguished vectors (e.g.\
singular vectors for the Virasoro algebra in the Fock space) on highest weight modules.
The key method is based on a~simple observation that the zero mode operator
of the vector $\exp\big({-}\sum\limits_{n>0}\frac{h(-n)}{n}z^{n}\big){\bf 1}$, which lives in the completed Fock space, acts on the
highest weight Fock module vector $v_\lambda$ semisimply with eigenvalue $(1+z)^{\langle h, \lambda \rangle}$.

\item Extra terms coming from Zhu's algebra multiplication: $a * b=\operatorname{Res}_{x_0} \frac{(1+x_0)^{{\deg}(a)}}{x_0}
Y(a,x_0)b$, namely $(1+x_0)^{{\rm deg}(a)}$.
\end{itemize}

Combination of these three methods for purposes of proving $C_2$-cof\/initeness and description of Zhu's algebra has led
to many constant term evaluations and identities.
Already a~sample of new identities can be obtained from the triplet vertex algebra~\cite{am1,am-sigma, am-imrn,
am-ja,am2}.
It is interesting that classif\/ication of modules for the triplet algebra can be reduced to a~single constant term
identity~\cite{am1}.

In this paper, we continue our investigation of the $ADE$ subalgebras of $\mathcal{W}(p)$ initiated in~\cite{ALM} in
connection to constant term identities.

We begin with a~short review of the triplet algebra $\mathcal{W}(p)$, its orbifold subalgebras $\mathcal{W}(p)^\Gamma$,
and status of classif\/ications of their modules; more details can be found in~\cite{ALM,ALM1,am1}.

Let $L=\mathbb{Z}\alpha$ be a~rank one lattice with $\langle\alpha,\alpha\rangle=2p$ ($p\geq1$).
Let
\begin{gather*}
V_{L}=\mathcal{U}(\widehat{\mathfrak{h}}_{<0})\otimes\mathbb{C}[{L}]
\end{gather*}
be the corresponding lattice vertex operator algebra~\cite{ll}, where $\mathfrak{\widehat{h}}$ is the af\/f\/inization of
$\mathfrak{h}=\mathbb{C}\alpha$, and let $\mathbb{C}[{L}]$ be the group algebra of~$L$.
Let $M(1)$ be the Heisenberg vertex subalgebra of $V_L$ generated by $\alpha(-1) \textbf{1}$.
With conformal vector
\begin{gather*}
\omega=\frac{\alpha(-1)^{2}}{4p}\textbf{1}+\frac{p-1}{2p}\alpha(-2)\textbf{1},
\end{gather*}
the space $V_{L}$ has a~vertex operator algebra structure of central charge
\begin{gather*}
c_{p}=1-\frac{6(p-1)^{2}}{p},
\end{gather*}
for more details see~\cite{ll}.
This vertex algebra admits two (degree one) screenings: $Q=e^{\alpha}_0$ and $\widetilde{Q}=e^{-\alpha/p}_0$, where we
used $e^\beta$ to denote vectors in the group algebra of the dual lattice of~$L$.
The triplet vertex algebra $\mathcal{W}(p)$ is def\/ined to be the kernel of the ``short'' screening $\widetilde{Q}$ on $V_L$.
It is strongly generated by~$\omega$ and three primary vectors
\begin{gather*}
F=e^{-\alpha}, \qquad  H=Q F, \qquad  E=Q^2 F
\end{gather*}
of conformal weight $2 p-1$.

We use $\overline{M(1)}$ to denote the singlet vertex algebra (cf.~\cite{A-2003, am-1}).
It is realized as a~vertex subalgebra of $\mathcal{W}(p)$ generated by~$\omega$ and~$H$.
Set
\begin{gather*}
e=Q,
\qquad
h=\frac{\alpha(0)}{p},
\end{gather*}
acting on $\mathcal{W}(p)$.
Let $f\in \End_{\Vir}(\mathcal {W}(p))$ be the unique operator def\/ined by
\begin{gather*}
fe^{-n\alpha}=0,
\qquad
fQ^{i}e^{-n\alpha}=i(2n+1-i)Q^{i-1}e^{-n\alpha},
\qquad
1\leq i\leq2n.
\end{gather*}
It was proved f\/irst in~\cite{ALM} that $\mathcal {W}(p)$ admits an action of $\mathfrak{sl}_{2}$ by the above three
operators.
The integration of the action of $\mathfrak{sl}_{2}$ gives rise to an action of ${\rm PSL}(2,\mathbb{C})$ on the vertex
operator al\-geb\-ra~$\mathcal {W}(p)$.
Vertex algebra $\mathcal{W}(p)^{A_m}$ (resp.\
$\mathcal{W}(p)^{D_m}$) is def\/ined as the invariant subalgebras with respect to the cyclic group $A_m$ of order~$m$
(reps.
dihedral group~$D_m$ of order~$2m$).
It is important to notice that $\mathcal{W}(p)^{D_m}$ is the $\mathbb{Z}_2$-orbifold of $\mathcal{W}(p)^{A_m}$ with
respect to the automorphism~$\Psi$ which is uniquely determined by the property
\begin{gather*}
\Psi \big(Q^i e^{- n \alpha}\big)=\frac{(-1)^i i!}{(2n-i)!} Q^{2n -i} e^ {-n \alpha}.
\end{gather*}
Note also that~$\Psi$ is also an automorphism of $\overline{M(1)}$ such that $\Psi (H)=-H$ and its f\/ixed point
subalgebra $\overline{M(1)}^+$ is a~subalgebra of $ \mathcal {W}(p)^{D_m}$ (for details see~\cite{ALM, ALM1}).

In~\cite{ALM,ALM1}, based on investigation of modules and twisted modules of $\mathcal{W}(p)$, we gave a~conjectural
list of $2m^2p$ irreducible $\mathcal {W}(p)^{A_m}$-modules and a~list of $(m^2+7)p$ irreducible $\mathcal
{W}(p)^{D_m}$-modules.
Moreover, for $m=2$, we showed that these two lists are complete under the assumption of validity of certain constant
term identities which have been verif\/ied by computer for small value~$p$.

In this paper we f\/irst give a~detailed investigation of Zhu's algebras of $\mathcal{W}(p)^{D_m}$ in connection to
classif\/ication of irreducible $\mathcal{W}(p)^{D_m}$-modules.
Results from~\cite{ALM1} show that Zhu's algebra $A(\mathcal{W}(p)^{D_m})$ is a~commutative, f\/inite-dimensional algebra
such that $\dim A(\mathcal{W}(p)^{D_m}) \ge (m^2 + 8) p -1$.
In order to prove that modules constructed in~\cite{ALM1} provide a~complete list of irreducible
$\mathcal{W}(p)^{D_m}$-modules, it suf\/f\/ices to prove inequality
\begin{gather*}
\dim A\big(\mathcal{W}(p)^{D_m}\big) \le \big(m^2 + 8\big) p -1.
\end{gather*}
We apply the methods similar to those used in~\cite{ALM} (see also~\cite{A-2003, am1, am2}) and evaluate certain
relations in Zhu's algebras on modules for Heisenberg vertex algebra $M(1)$.
This leads to a~new series of constant term identities which we list in Appendix~\ref{ct-section}.

The structure of Zhu's algebra $A(\mathcal{W}(p)^{A_m}) $ is discussed in Section~\ref{class-A}.
We f\/irst show that the classif\/ication of irreducible modules for $\mathcal{W}(p)^{A_m} $ is equivalent to proving that
Zhu's algebra $A(\mathcal{W}(p)^{A_m}) $ has dimension $(2 m^2 + 4) p -1$.
The results from~\cite{ALM} implies that $\dim A(\mathcal{W}(p)^{A_m}) \ge (2 m^2 + 4) p -1$.
In Section~\ref{class-A} we present a~detailed discussion of the proof of the opposite inequality, and f\/inally show that
it can be reduced to the proof of certain combinatorial identities which we checked for small values of~$m$ and~$p$.

Throughout this paper we use
\begin{gather*}
h_{i,j}=\frac{(i-jp+p-1)(i-jp-p+1)}{4p}
\end{gather*}
to parameterize lowest conformal weights of modules.

\section[On classif\/ication of irreducible $\mathcal{W}(p)^{D_m}$-modules]{On classif\/ication
of irreducible $\boldsymbol{\mathcal{W}(p)^{D_m}}$-modules}\label{Section2}

In~\cite{ALM1}, we initiated the study of representation theory of the vertex operator algebra $\mathcal{W}(p)^{D_m}$,
$m \geq 2$.
We proved $C_2$-cof\/initeness of these vertex algebras and showed that the associated Zhu's algebra is
a~f\/inite-dimensional commutative algebra.
In the representation theory of $\mathcal{W}(p)^{D_m}$, the singlet vertex algebra $\overline{M(1)}^{+}$ has played an
important role.
In the same paper, we classif\/ied all irreducible modules for $\overline{M(1)}^{+}$ and $\mathcal{W}(p)^{D_2}$ modulo the
same constant term identity.

In this section we shall slightly extend results from~\cite{ALM1} so we shall introduce two new combinatorial
conjectures which will imply the classif\/ication of irreducible modules.

First we recall some results we obtained in~\cite{ALM1}.
\begin{Theorem}\quad
\begin{enumerate}\itemsep=0pt
\item[$(1)$]The vertex algebra $\mathcal{W}(p)^{D_m}$ is strongly generated~by
\begin{gather*}
\omega,
\qquad
H^{(2)}=Q^{2}e^{-2\alpha},
\qquad
U^{(m)}=(2m)! F^{(m)} + E^{(m)},
\end{gather*}
where $F^{(m)}=e^{- m \alpha}$, $ E^{(m)}=Q^{2 m} F^{(m)}$.
\item[$(2)$] Zhu's algebra $A(\mathcal{W}(p)^{D_m})$ is a~commutative associative algebra generated by $[\omega]$,
$[H^{(2)}]$ and $[U^{(m)}]$, where $[\enskip \cdot \enskip]$ denotes coset of an element in Zhu's algebra.

\item[$(3)$] $\mathcal{W}(p)^{D_m}$ has $(m^2 + 7) p$ inequivalent irreducible modules constructed from twisted and
untwisted $\mathcal{W}(p)$-modules whose lowest components are all $1$-dimensional.
\end{enumerate}
\end{Theorem}

Since $\mathcal{W}(p)^{D_m}$ has $p-1$ logarithmic modules constructed in~\cite{am2}, we conclude that
$A(\mathcal{W}(p)^{D_m})$ also has $p-1$ indecomposable $2$-dimensional modules.
So we have:

\begin{Corollary}
For every $m \ge 2$
\begin{gather*}
\dim A\big(\mathcal{W}(p)^{D_m}\big) \ge \big(m^2 + 8\big) p-1.
\end{gather*}
\end{Corollary}

In order to classify irreducible modules it is suf\/f\/icient to prove the inequality
\begin{gather}
\dim A\big(\mathcal{W}(p)^{D_m}\big) \le \big(m^2 + 8\big) p-1.
\label{nejednakost}
\end{gather}
This will imply that the dimension of Zhu's algebra is of course $(m^2 + 8) p -1$ and therefore the list of irreducible
modules constructed in~\cite{ALM1} (see also Tables 1 and 2) will be a~complete list of irreducible
$\mathcal{W}(p)^{D_m}$-modules, up to equivalence.
We obtained this result in~\cite{ALM1} for $m=2$.

\begin{Theorem}
Assume that Conjecture {\rm 7.6} of~{\rm \cite{ALM1}} holds $($verified by computer for small~$p)$.
Then
\begin{gather*}
\dim A\big(\mathcal{W}(p)^{D_2}\big)=12 p-1,
\end{gather*}
and $\mathcal{W}(p)^{D_2}$ has exactly $11 p$ irreducible modules constructed explicitly in~{\rm \cite{ALM1}}.
\end{Theorem}

Next we shall see that for general~$m$ the inequality~\eqref{nejednakost} is also related to certain constant term
combinatorial identities which are more complicated than the identities which appear in~\cite{ALM1}.

Note that $\mathcal{W}(p)^{D_m}$ contains a~subalgebra generated by~$\omega$ and $H^{(2)}$.
As in~\cite{ALM1}, we denoted this subalgebra by $\overline{M(1)}^+$ since it is a~$\mathbb{Z}_2$ orbifold of the
singlet vertex algebra $\overline{M(1)}$~\cite{A-2003}, which is generated by~$\omega$ and $H=Q e^{-\alpha}$.
We also determined the structure of Zhu's algebra $A (\overline{M(1)}^+) $.
The next result is again taken from~\cite{ALM1}.

\begin{Proposition}\qquad
\begin{enumerate}\itemsep=0pt
\item[$(1)$] Inside the Zhu algebra $A (\overline{M(1)}^+) $ we have the following relations:
\begin{gather*}
\big(\big[H^{(2)}\big]-f_p([\omega])\big) * \big(\big[H^{(2)}\big]-r_p([\omega])\big)=0,
\end{gather*}
where $r_p \in {\mathbb C}[x]$, $\deg r_p \le 3p-1$, and
\begin{gather*}
f_p(x)=(-1)^{p}\frac{(4p)^{3p-1}(2p)!}{(4p-1)!(3p-1)!p!} \prod\limits_{i=1}^{3p-1} (x- h_{i,1}).
\end{gather*}

\item[$(2)$] Assume that Conjecture~{\rm \ref{conjecutreD2}} holds, then in $A (\overline{M(1)}^+) $
\begin{gather*}
\ell_p([\omega]) * \big(\big[H^{(2)}\big]-f_p([\omega])\big)=0,
\end{gather*}
where
\begin{gather*}
\ell_p(x)=\prod\limits_{i=1}^{p}(x- h_{4p-i,1}) \prod\limits_{i=1}^{2p} (x-h_{3p+1/2 -i,1}).
\end{gather*}
\end{enumerate}
\end{Proposition}

{\sloppy Let $A_{0}(\mathcal{W}(p)^{D_m})$ be the subalgebra of $A(\mathcal{W}(p)^{D_m})$ generated by $[\omega]$ and
$[H^{(2)}]$.
Let $A_{1}(\mathcal{W}(p)^{D_m}) =A_{0}(\mathcal{W}(p)^{D_m}).
[U^{(m)}]$.

}

\begin{Lemma}
\label{uu}
In $A(\mathcal{W}(p)^{D_m})$, we have $[U^{(m)}]*[U^{(m)}]\in A_{0}(\mathcal{W}(p)^{D_m})$.
Moreover, $A(\mathcal{W}(p)^{D_m})$ is a~$\mathbb{Z}_2$-graded algebra
\begin{gather*}
A\big(\mathcal{W}(p)^{D_m}\big)=A_0\big(\mathcal{W}(p)^{D_m}\big) \oplus A_1\big(\mathcal{W}(p)^{D_m}\big).
\end{gather*}
\end{Lemma}
\begin{proof}
First we notice that
\begin{gather*}
E^{(m)} * E^{(m)}=F^{(m)} * F^{(m)}=0,
\qquad
E^{(m)} * F^{(m)},
\qquad
F^{(m)} * E^{(m)} \in \overline{M(1)},
\end{gather*}
which implies that $U^{(m)}* U^{(m)} \in \overline{M(1)}^{+}$.
The proof follows.
\end{proof}

For technical reasons we need to recall some informations on lowest weights of irreducible modules.
Since Zhu's algebra $A(\mathcal{W}(p)^{D_m})$ is commutative (see below), its irreducible modules are $1$-dimensional.
Then applying Zhu's algebra theory we see that all irreducible $A(\mathcal{W}(p)^{D_m})$-modules should be
parameterized by its lowest weights with respect to $(L(0), H^{(2)} (0), U^{(m)} (0))$.
The lowest weights of irreducible $\mathcal{W}(p)^{D_m}$-module constructed in~\cite{ALM1} can be found in the
following two tables, where
\begin{gather*}
\phi(t)=(-1)^{\frac{m(m-1)p}{2}}\prod\limits_{l=0}^{m-1}{t+pl \choose
(m+1)p-1}\frac{((m+1)p-1)!((l+1)p)!}{((m+l+1)p-1)!p!},
\\
\sigma=\big[\tfrac{1}{2}(
\begin{smallmatrix}
1+i & -1+i
\\
1+i & 1-i
\end{smallmatrix}
)\big]\in {\rm PSL}(2,\mathbb{C}),
\end{gather*}
and the number~$\ell$ is def\/ined as in~\cite[Lemma~4.8]{ALM}.

\begin{table}[t] \centering \small \caption{Irreducible $\mathcal{W}(p)^{D_m}$-modules: $m=2k$.}\label{table1}\vspace{1mm}
\begin{tabular}{|c|c|c|c|}
\hline
\mbox{module}~$M$ & $L(0)$ & $ H^{(2)}(0) $ & $ U^{(m)}(0) $
\\
\hline
$\Lambda(i)^+_0 $\tsep{1pt}\bsep{1pt} & $ h_{i,1}$& $ 0 $ & $0 $
\\
\hline
$\Lambda(i)^-_0 $\tsep{1pt}\bsep{1pt} & $ h_{i,3}$& $ -2f_p (h_{i,3}) $ & $0 $
\\
\hline
$ \Lambda(i)_j^ +\cong\Lambda(i)_j^ - $\tsep{1pt}\bsep{1pt} & $ h_{i, 2 j +1} $& $ f_p (h_{i, 2 j +1}) $ & $0 $
\\
\hline
$ \Lambda(i)^+_m $\tsep{2pt}\bsep{2pt} & $ h_{i, m +1} $ & $ f_p (h_{i, m+1}) $ & $\frac{(2m)!}{m!}\phi(-mp+i-1) $
\\
\hline
$ \Lambda(i)^-_m $\tsep{2pt}\bsep{2pt} & $ h_{i, m +1} $ & $ f_p (h_{i, m+1}) $ & $-\frac{(2m)!}{m!}\phi(-mp+i-1) $
\\
\hline
$ \Pi(i)_j^+ \cong\Pi(i)_j^- $\tsep{1pt}\bsep{1pt} & $ h_{p+i,2 j+1} $ & $ f_p (h_{p+i,2 j+1}) $ & $0$
\\
\hline
$R(i, j,k)$\tsep{1pt}\bsep{1pt} & $h_{\ell + 1-i/m,1}$ &$f_p (h_{\ell + 1-i/m,1})$& $0$
\\
\hline
$R(j)^{\sigma}$\tsep{2pt}\bsep{2pt} & $ h_{3p + 1/2- j,1}$ &$- \frac{1}{2} f_p (h_{3p+1/2-j,1})$& $\frac{i^m(2m)!}{2^{m-1}m!}\phi(3p-1/2-
j)$
\\
\hline
$R(j)^{h \sigma}$\tsep{2pt}\bsep{2pt} & $ h_{3p + 1/2-j,1}$ &$-\frac{1}{2} f_p (h_{3p+1/2-j,1}) $& $-\frac{i^m(2m)!}{2^{m-1}m!}\phi(3p-1/2- j)$
\\
\hline
\end{tabular}
\end{table}

\begin{table}[t] \small \centering\caption{Irreducible $\mathcal{W}(p)^{D_m}$-modules: $m=2k+1$.}\label{table2}\vspace{1mm}
\begin{tabular}{|c|c|c|c|}
\hline
\mbox{module}~$M$ & $L(0)$ & $ H^{(2)}(0) $ & $ U^{(m)}(0) $
\\
\hline
$\Lambda(i)^+_0 $\tsep{1pt}\bsep{1pt}  & $ h_{i,1}$& $ 0 $ & $0 $
\\
\hline
$\Lambda(i)^-_0 $\tsep{1pt}\bsep{1pt}  & $ h_{i,3}$& $ -2 f_p (h_{i,3}) $ & $0 $
\\
\hline
$ \Lambda(i)_j^ +\cong\Lambda(i)_j^ - $\tsep{1pt}\bsep{1pt}  & $ h_{i, 2 j +1} $& $ f_p (h_{i, 2 j +1}) $ & $0 $
\\
\hline
$ \Pi(i)_j^+ \cong\Pi(i)_j^- $\tsep{1pt}\bsep{1pt}  & $ h_{p+i,2 j+1} $ & $ f_p (h_{p+i,2 j+1}) $ & $0$
\\
\hline
$ \Pi(i)^+_m $\tsep{2pt}\bsep{2pt}  & $ h_{p+i, m+2} $ & $ f_p (h_{p+i, m+2}) $ & $\frac{(2m)!}{m!}\phi(-mp+i-1) $
\\
\hline
$ \Pi(i)^-_m $\tsep{2pt}\bsep{2pt}  & $ h_{p+i, m+2} $ & $ f_p (h_{p+i, m+2})$ & $-\frac{(2m)!}{m!}\phi(-mp+i-1)$
\\
\hline
$R(i, j,k)$\tsep{1pt}\bsep{1pt}  & $h_{\ell + 1-i/m,1}$&$ f_p (h_{\ell + 1-i/m,1})$& $0$
\\
\hline
$R(j)^{\sigma}$\tsep{2pt}\bsep{2pt}  & $ h_{3p + 1/2- j,1} $&$- \frac{1}{2} f_p (h_{3p+1/2-j,1})$& $\frac{i^m(2m)!}{2^{m-1}m!}\phi(3p-1/2-
j)$
\\
\hline
$R(j)^{h \sigma}$\tsep{2pt}\bsep{2pt}  & $ h_{3p + 1/2-j,1}$&$-\frac{1}{2} f_p (h_{3p+1/2-j,1}) $& $-\frac{i^m(2m)!}{2^{m-1}m!}\phi(3p-1/2- j)$
\\
\hline
\end{tabular}
\end{table}

By the above two tables and standard arguments we infer:

\begin{Lemma}
We have the following relation in $A(\mathcal{W}(p)^{D_m})$:
\begin{gather*}
\big[H^{(2)}\big]*\big[U^{(m)}\big]=\big[U^{(m)}\big]*\big[H^{(2)}\big]=h_p([\omega])\big[U^{(m)}\big],
\end{gather*}
where $h_p(x)$ is a~polynomial of degree at most $3p-1$, and satisfies the following interpolation conditions:
\begin{gather*}
h_p(h_{i,m+1})=f_p (h_{i,m+1}),
\qquad
h_{p} (h_{3p+1/2-j,1})=- \tfrac{1}{2} f_p (h_{3p+1/2-j,1}),
\end{gather*}
for $i=1, \dots, p$, $j=1, \dots, 2p$.
In particular, Zhu's algebra $A(\mathcal{W}(p)^{D_m})$ is commutative.
\end{Lemma}
\begin{proof}
We use the above tables and apply both vectors in the equation on lowest weight vectors of modules.
\end{proof}

Next, for $a, b \in \mathcal{W}(p)$ we def\/ine
\begin{gather*}
a{\tilde{\circ}} b=\Res_z \frac{(1+z)^{\deg (a)}} {z^3} Y(a,z) b
\qquad \text{and}\qquad
a{\tilde{\circ}}_k b=\Res_z \frac{(1+z)^{\deg (a)}} {z^k} Y(a,z) b.
\end{gather*}

\begin{Lemma}
\label{HU1}
Assume that Conjecture~{\rm \ref{conjecutreDm}} holds for $m \ge 3$ and that Conjecture~{\rm \ref{conjecutreD2}} holds for $m=2$.
Then in $A(\mathcal{W}(p)^{D_m})$,
\begin{gather*}
g_p([\omega])*\big[U^{(m)}\big]=0,
\end{gather*}
where
\begin{gather*}
g_p(x)=\prod\limits_{i=1}^{p}(x- h_{i,m+1}) \prod\limits_{i=1}^{2p} (x-h_{3p+1/2 -i,1}).
\end{gather*}
In particular,
\begin{gather*}
\dim A_1 \big(\mathcal{W}(p)^{D_m}\big) \le 3p.
\end{gather*}
\end{Lemma}

\begin{proof}
We assume that $m\geq3$, noticing that the case $m=2$ has been treated in~\cite{ALM1}.
From the structure of $\mathcal{W}(p)$ as a~Virasoro module
\begin{gather*}
U^{(m)} {\tilde{\circ}} H^{(2)}=t_{p}U^{(m)},
\end{gather*}
where $t_{p}\in\mathcal{U}(Vir^{-})$.
Then
\begin{gather*}
Q^{m}\big(e^{-m\alpha} {\tilde{\circ}} H^{(2)}\big)=t_{p}Q^{m}e^{-m\alpha}.
\end{gather*}
We see that there exists a~polynomial $ g_p(x)$ of degree at most~$3p$, such that
\begin{gather*}
g_p([\omega])*\big[U^{(m)}\big]=0
\end{gather*}
in $A(\mathcal{W}(p)^{D_m})$, and
\begin{gather*}
\big[Q^{m}\big(e^{-m\alpha} {\tilde{\circ}} H^{(2)}\big)\big]=g_{p}([\omega])*[Q^{m}e^{-m\alpha}]
\end{gather*}
in $A(\overline{M(1)})$.
By evaluating the left hand side on known irreducible $\mathcal{W}(p)^{D_m}$-modules we see that
\begin{gather*}
g_p(x)=C_p \prod\limits_{i=1}^{p}(x- h_{i,m+1}) \prod\limits_{i=1}^{2p} (x-h_{3p+1/2 -i,1})
\end{gather*}
for some constant $C_p $.

As in our previous papers we shall relate evaluation of the constant $C_p$ with action of certain elements of
$\overline{M(1)}$ on lowest weight $\overline{M(1)}$-modules.
These highest weight modules are realized as modules $M(1,\lambda)$ for the Heisenberg vertex algebra $M(1)$ (remember
that $\overline{M(1)} \subset M(1)$).

So let $v_{\lambda}$ be a~lowest weight vector in the $M(1, \lambda)$ and let $t=\langle \lambda, \alpha \rangle $.
Then we get
\begin{gather*}
o \big(Q^{m}\big(e^{-m\alpha} {\tilde{\circ}} H^{(2)}\big)\big) v_{\lambda}
\\
\qquad
={\rm Res}_{x_{0},x_{1},\dots, x_{m+2}}\frac{(1+x_{0})^{m^2p+mp-(t+1)m}} {x_{0}^{-4mp+3}(x_{1}\cdots x_{m+2})^{4p}}(1+x_{1})^{t}\cdots(1+ x_{m+2})^{t}
\\
\qquad
\phantom{=}{}
\times
(x_{0}-x_{m+1})^{-2mp}(x_{0}-x_{m+2})^{-2mp}\prod\limits_{1\leq i<j\leq m+2}(x_{i}-x_{j})^{2p}\prod\limits_{i=1}^{m}(x_{i}-x_{0})^{-2mp}.
\end{gather*}
It follows from Conjecture~\ref{conjecutreDm} that $C_p $ is nonzero.
\end{proof}

\begin{Remark}
Notice that the previous lemma generalizes our result from~\cite{ALM1}.
There we proved that for $m=2$, $\ell_p([\omega])*[U^{(2)}]=0$.
Observe that $\ell_p=g_p$ only for $m=2$.
\end{Remark}

Now we want to calculate upper bound for $\dim A_1 (\mathcal{W}(p)^{D_m})$.

Using a~similar calculation as above, we have
\begin{gather*}
\big[U^{(m)}{\tilde{\circ}} U^{(m)}\big]=k_p([\omega])\big(\big[H^{(2)}\big]-f_p([\omega])\big)
\\
\phantom{\big[U^{(m)}{\tilde{\circ}} U^{(m)}\big]=}{}
+l_p([\omega]) \prod\limits_{i=1}^{p}([\omega]- h_{i,1})\prod\limits_{i=1}^{p-1}([\omega]- h_{mp+p+i,1})
\prod\limits_{i=1}^{m^2p} \bigg([\omega]-h_{p+\tfrac{i}{m},1}\bigg),
\end{gather*}
and
\begin{gather*}
\big[U^{(m)}{\tilde{\circ}}_{5} U^{(m)}\big]=\tilde{k}_p([\omega])\big(\big[H^{(2)}\big]-f_p([\omega])\big)
\\
\phantom{\big[U^{(m)}{\tilde{\circ}}_{5} U^{(m)}\big]=}{}
+\tilde{l}_p([\omega])  \prod\limits_{i=1}^{p}([\omega]- h_{i,1})\prod\limits_{i=1}^{p-1}\big([\omega]- h_{mp+p+i,1}\big)
\prod\limits_{i=1}^{m^2p} \bigg([\omega]-h_{p+\tfrac{i}{m},1}\bigg),
\end{gather*}
where $k_p(x)$, $\tilde{k}_p(x)$, $l_p(x)$, $\tilde{l}_p(x)\in {\mathbb C}[x]$, and
\begin{gather*}
\deg l_p \le (m-2)(p-1),
\qquad
\deg \tilde{l}_p \le (m-2)(p-1)+1.
\end{gather*}

By the proof of~\cite[Lemma 2.1.3]{z}, we get
\begin{gather*}
U^{(m)}{\tilde{\circ}} U^{(m)}=\Res_x \frac{(1+x)^{m^2p+mp-m} (2+x)} {x^3} Y(e^{-m\alpha},x) Q^{2m}e^{-m\alpha},
\end{gather*}
and
\begin{gather*}
U^{(m)}{\tilde{\circ}}_{5} U^{(m)}=\Res_x \frac{(1+x)^{m^2p+mp-m} (2+3x+3x^{2}+x^{3})} {x^5} Y(e^{-m\alpha},x)Q^{2m}e^{-m\alpha}.
\end{gather*}

Hence
\begin{gather*}
o\big(U^{(m)}{\tilde{\circ}} U^{(m)}\big) v_{\lambda}=H_{p,m}(t) v_{\lambda}
\qquad \text{and}\qquad
o\big(U^{(m)}{\tilde{\circ}}_{5} U^{(m)}\big) v_{\lambda}=\tilde{H}_{p,m}(t) v_{\lambda},
\end{gather*}
where
\begin{gather*}
H_{p,m}(t)={\rm Res}_{x_{0},\dots,x_{2m}}\frac{(1+x_{0})^{m^2p+mp-(t+1)m} (2+x_{0})} {x_{0}^{2m^{2}p+3}}(x_{1}\cdots
x_{2m})^{-2mp}
\\
\phantom{H_{p,m}(t)=}{}
\times\prod\limits_{i=1}^{2m} (1+ x_{i})^{t} \prod\limits_{1\leq i<j\leq 2m}(x_{i}-x_{j})^{2p}\prod\limits_{i=1}^{2m}
\left(1-\frac{x_{i}}{x_{0}}\right)^{-2mp},
\end{gather*}
and
\begin{gather*}
\tilde{H}_{p,m}(t)={\rm Res}_{x_{0}, \cdot, x_{2m}}\frac{(1+x_{0})^{m^2p+mp-(m+1)t} \big(2+3x_{0}+3x_{0}^{2}+x_{0}^{3}\big)}
{x_{0}^{2m^{2}p+5}}(x_{1}\cdots x_{2m})^{-2mp}\\
\phantom{\tilde{H}_{p,m}(t)=}{}
\times\prod\limits_{i=1}^{2m} (1+ x_{i})^{t} \prod\limits_{1\leq i<j\leq 2m}(x_{i}-x_{j})^{2p}\prod\limits_{i=1}^{2m}
\left(1-\frac{x_{i}}{x_{0}}\right)^{-2mp}.
\end{gather*}

We can show that for small values of $(m,p)$
\begin{gather*}
H_{p,m}(t)=h_{p,m}(t){t \choose p}{t+1-p \choose p}{t-(m+1)p \choose p -1}{t+mp \choose p -1}
\prod\limits_{i=1}^{m^2p} \left((t+1-p)^2-\frac{i^2}{m^2}\right),\\
\tilde{H}_{p,m}(t)=\tilde{h}_{p,m}(t){t \choose p}{t+1-p \choose p}{t-(m+1)p \choose p -1}{t+mp \choose p -1}
\prod\limits_{i=1}^{m^2p} \left((t+1-p)^2-\frac{i^2}{m^2}\right),
\end{gather*}
where $h_{p,m}$ and $\tilde{h}_{p,m}$ are given in the table (up to a~scalar factor):

\begin{center}
\begin{tabular}{|c|c|c|c}
\hline
$(m, p)$ & $h_{p,m}(t)$ & $ \tilde{h}_{p,m}(t) $\tsep{1pt}
\\
\hline
 $(2, 1) $& $ 1 $ & $4 t^2+107 $\tsep{1pt}
\\
\hline
$(2, 2) $& $ 1 $ & $4 t^2-8 t+175 $\tsep{1pt}
\\
\hline
$ (2, 3) $ & $ 1 $ & $4 t^2-16 t+219 $\tsep{1pt}
\\
\hline
$ (2, 4) $ & $ 1 $ & $4 t^2-24 t+239$\tsep{1pt}
\\
\hline
$(3, 1) $& $ 1 $ & $9 t^2+362 $\tsep{1pt}
\\
\hline
$(3, 2) $& $ t^2-2 t-30 $ & $549 t^4-2196 t^3+57020 t^2-109648 t-2118976$\tsep{1pt}
\\
\hline
$(3, 3) $& $ 26141 t^4-104564 t^3-576380 t^2 $ & $4977 t^6-29862 t^5+495920 t^4-1784600 t^3$\tsep{1pt}
\\
& ${}+1361888 t-3720960$ & ${}-13382432 t^2+30254432 t-84341760 $
\\
\hline
\end{tabular}
\end{center}

\begin{Conjecture}
\label{D1m}
Polynomials $h_{p,m}$ and $\tilde{h}_{p,m}$ are relatively prime.
\end{Conjecture}

This yields the following result:
\begin{Lemma}
Assume that Conjecture~{\rm \ref{D1m}} holds.
Then in $A(\mathcal{W}(p)^{D_m})$, we have
\begin{gather*}
k_p([\omega])\big(\big[H^{(2)}\big]-f_p([\omega])\big)+\prod\limits_{i=1}^{p}([\omega]- h_{i,1})\prod\limits_{i=1}^{p-1}([\omega]-
h_{mp+p+i,1}) \prod\limits_{i=1}^{m^2p} \bigg([\omega]-h_{p+\tfrac{i}{m},1}\bigg)=0
\end{gather*}
for some $k_p(x)\in {\mathbb C}[x]$.
\end{Lemma}

This lemma and the structure of Zhu's algebra $A(\overline{M(1)}^+)$ imply:

\begin{Proposition}
\label{dim}
Assume that Conjectures~{\rm \ref{D1m}} and~{\rm \ref{conjecutreDm}} hold.
Then
\begin{gather*}
\dim A_{0}\big(\mathcal{W}(p)^{D_m}\big) \le m^2p+5p-1.
\end{gather*}
\end{Proposition}

Now we are in a~position to give the f\/irst main result of this paper.
\begin{Theorem}
Assume that Conjectures~{\rm \ref{D1m}},~{\rm \ref{conjecutreDm}} and~{\rm \ref{conjecutreD2}} hold.
Then
\begin{enumerate}\itemsep=0pt
\item[$(1)$] Tables~{\rm \ref{table1}} and~{\rm \ref{table2}} give a~complete list of irreducible $\mathcal{W}(p)^{D_m}$-modules.
In particular, $\mathcal{W}(p)^{D_m}$ has $(m^2 + 7) p$ non-isomorphic irreducible modules,

\item[$(2)$] Zhu's algebra $A(\mathcal{W}(p)^{D_m})$ is a~commutative algebra of dimension $(m^2 + 8)p-1$.
\end{enumerate}
\end{Theorem}

\begin{proof}
As we discussed above it is enough to prove inequality~\eqref{nejednakost}.
We proved in~\cite{ALM1} that $A(\mathcal{W}(p)^{D_m})$ is commutative.
We will compute the dimension of $A(\mathcal{W}(p)^{D_m})$.
By Lemma~\ref{HU1} and Proposition~\ref{dim}, we get that
\begin{gather*}
\dim A\big(\mathcal{W}(p)^{D_m}\big)=\dim A_{0}\big(\mathcal{W}(p)^{D_m}\big) + \dim A_1\big(\mathcal{W}(p)^{D_m}\big)
\\
\phantom{\dim A\big(\mathcal{W}(p)^{D_m}\big)}
 \le  m^2 p + 5 p -1 + 3p=\big(m^2 + 8\big) p -1.
\end{gather*}
The proof follows.
\end{proof}

\section[On classif\/ication of irreducible $\mathcal{W}(p)^{A_m}$-modules]{On classif\/ication
of irreducible $\boldsymbol{\mathcal{W}(p)^{A_m}}$-modules}\label{class-A}

In our paper~\cite{ALM} we constructed $2 m^2 p$ irreducible $\mathcal{W}(p)^{A_m}$-modules and conjectured that these
modules provide a~complete list of irreducible $\mathcal{W}(p)^{A_m}$-modules.
In the case $m=2$ we presented a~proof which is based on certain constant term identities.
These identities are dif\/f\/icult to prove in general, but using Mathematica they can be verif\/ied for~$p$ small.
So our approach can be considered as an algorithm that reduces problems in representation theory to checking something
purely computational.

In this section we shall extend our results of~\cite{ALM} to general~$m$ and thus provide an algorithm for
a~classif\/ication of irreducible $\mathcal{W}(p)^{A_m}$-modules.

Let us f\/irst recall some results from~\cite{ALM}.
Set $F^{(m)}=e^{-m\alpha}$, $E^{(m)}=Q^{2m}e^{-m\alpha}$, $H=Qe^{-\alpha}$. Then $\mathcal {W}(p)^{A_{m}}$ is strongly
generated by $E^{(m)}$, $F^{(m)}$, $H$ and~$\omega$.
Hence, Zhu's algebra $A(\mathcal{W}(p)^{A_m})$ is generated by $[\omega]$, $[H]$, $[E^{(m)}]$ and~$[F^{(m)}]$.

It is also important to notice that the restriction of the automorphism~$\Psi$ of $\mathcal{W}(p)$ (cf.~\cite{ALM}) to
$\mathcal{W}(p)^{A_m}$ gives an automorphism of order two such that
\begin{gather*}
\Psi (F^{(m)})=E^{(m)},
\qquad
\Psi(H)=- H.
\end{gather*}

In~\cite{ALM} we constructed $ 2 m^2 p$ irreducible $\mathcal{W}(p)^{A_m}$-modules.
A~list of irreducible $\mathcal{W}(p)^{A_m}$-modules and their lowest weights with respect to $(L(0), H(0))$ are given
in the following tables:

\begin{table}[htp]\small \centering\caption{Irreducible $\mathcal{W}(p)^{A_m}$-modules: $m=2k$.}\label{table3}\vspace{1mm}
\begin{tabular}{|c|c|c|c}
\hline
\mbox{module}~$M$ & \mbox{lowest weights} & $ \dim M(0) $\tsep{1pt}\bsep{1pt}
\\
\hline
$\Lambda(i)_0 $ & $ (h_{i,1}, 0) $\tsep{1pt}\bsep{1pt} & $1 $
\\
\hline
$ \Lambda(i)_j^ + $ & $ \big(h_{i, 2 j +1}, {- 2 j p-1 + i \choose 2p-1}\big) $\tsep{3pt}\bsep{3pt} & $1$
\\
\hline
$ \Lambda(i)_j^ - $ & $ \big(h_{i, 2 j +1}, -{- 2 j p-1 + i \choose 2p-1}\big) $\tsep{3pt}\bsep{3pt} & $1$
\\
\hline
$ \Lambda(i)_m $ & $ \big(h_{i, 2 k +1}, \pm{- 2 k p-1 + i \choose 2p-1}\big) $\tsep{3pt}\bsep{3pt} & $2$
\\
\hline
$ \Pi(i)_j^+ $ & $ \big(h_{p+i,2 j+1}, {-(2j -1) p-1 + i \choose 2 p -1}\big) $\tsep{3pt}\bsep{3pt} & $1$
\\
\hline
$\Pi(i)_j^- $ & $ \big(h_{p+i,2 j +1},-{-(2j-1) p-1 + i \choose 2 p -1}\big) $\tsep{3pt}\bsep{3pt} & $1 $
\\
\hline
$R(i, j,k)$ & $\big(h_{\ell + 1-i/m,1}, {\ell-\tfrac{i}{m} \choose 2p-1}\big)$\tsep{3pt}\bsep{3pt}& $1$
\\
\hline
\end{tabular}
\end{table}

\begin{table}[htp]\small \centering\caption{irreducible $\mathcal{W}(p)^{A_m}$-modules: $m=2k+1$.}\vspace{1mm}
\begin{tabular}{|c|c|c|c}
\hline
\mbox{module}~$M$ & \mbox{lowest weights} & $ \dim M(0) $\tsep{1pt}\bsep{1pt}
\\
\hline
$\Lambda(i)_0 $ & $ (h_{i,1}, 0) $\tsep{1pt}\bsep{1pt} & $1 $
\\
\hline
$ \Lambda(i)_j^ + $ & $ \big(h_{i, 2 j +1}, {- 2 j p-1 + i \choose 2p-1}\big) $\tsep{3pt}\bsep{3pt} & $1$
\\
\hline
$ \Lambda(i)_j^ - $ & $ \big(h_{i, 2 j +1}, -{- 2 j p-1 + i \choose 2p-1}\big) $\tsep{3pt}\bsep{3pt} & $1$
\\
\hline
$ \Pi(i)_m $ & $ \big(h_{p+i, 2 k +3}, \pm{- (2 k+1) p-1 + i \choose 2p-1}\big) $\tsep{3pt}\bsep{3pt} & $2$
\\
\hline
$ \Pi(i)_j^+ $ & $ \big(h_{p+i,2 j+1}, {-(2j -1) p-1 + i \choose 2 p -1}\big) $\tsep{3pt}\bsep{3pt} & $1$
\\
\hline
$\Pi(i)_j^- $ & $ \big(h_{p+i,2 j+1},-{-(2j-1) p-1 + i \choose 2 p -1}\big) $\tsep{3pt}\bsep{3pt}& $1 $
\\
\hline
$R(i, j,k)$ & $ \big(h_{\ell + 1-i/m,1}, {\ell-\tfrac{i}{m} \choose 2p-1}\big)$\tsep{3pt}\bsep{3pt} & $1$
\\
\hline
\end{tabular}
\end{table}

\begin{Remark}
When $m=2k$, the action of $A(\mathcal{W}(p)^{A_m})$ on the lowest weight space of $\Lambda(i)_m$ is given~by
\begin{gather*}
  H(0)e^{- m\alpha/2+ \frac{(i-1) \alpha}{2p}}={- m p-1 + i \choose 2p-1}e^{- m\alpha/2+ \frac{(i-1) \alpha}{2p}},
\\
  H(0)Q^me^{- m\alpha/2+ \frac{(i-1) \alpha}{2p}}=-{- m p-1 + i \choose 2p-1}Q^me^{- m\alpha/2+ \frac{(i-1)
\alpha}{2p}},
\\
  E^{(m)}(0)e^{- m\alpha/2+ \frac{(i-1) \alpha}{2p}}=a_{m,p}Q^me^{- m\alpha/2+ \frac{(i-1) \alpha}{2p}},
\\
  F^{(m)}(0)Q^me^{- m\alpha/2+ \frac{(i-1) \alpha}{2p}}=b_{m,p}e^{- m\alpha/2+ \frac{(i-1) \alpha}{2p}},
\end{gather*}
where $a_{m,p}$, $b_{m,p}$ are nonzero constants.
Similar result holds for $\Pi(i)_m $ for $m=2k+1$.
\end{Remark}

So we constructed $(2 m^2-1) p$ irreducible modules whose lowest components are $1$-dimensional and~$p$ irreducible
modules with $2$-dimensional lowest components.
Since these lowest components are irreducible modules for Zhu's algebra $A(\mathcal{W}(p)^{A_m})$ and since
$A(\mathcal{W}(p)^{A_m})$ has $(p-1)$ $2$-di\-mensional
indecomposable modules constructed from logarithmic $\mathcal{W}(p)$-modules we get:
\begin{gather*}
\dim A\big(\mathcal{W}(p)^{A_m}\big) \ge \big(2 m^2 -1\big) p + 4p + p-1=\big(2 m^2 + 4\big) p-1.
\end{gather*}

So if we prove that in the above relation equality holds, we get the classif\/ication of irreducible modules.

\begin{Lemma}
\label{kriterij}
Assume that $\dim A(\mathcal{W}(p)^{A_m}) \le (2 m^2 + 4) p-1$.
Then the above two tables give a~complete list of irreducible $\mathcal{W}(p)^{A_m}$-modules.
\end{Lemma}

Let $A_{0}(\mathcal{W}(p)^{A_m})$ be the subalgebra of $A(\mathcal{W}(p)^{A_m})$ generated by $[\omega]$ and $[H]$, and let\linebreak
$A_{\omega}(\mathcal{W}(p)^{A_m})$ be the subalgebra generated by $[\omega]$.
Let
\begin{gather*}
A_{1}\big(\mathcal{W}(p)^{A_m}\big):=A_{0}\big(\mathcal{W}(p)^{A_m}\big).  E^{(m)},
\qquad
A_{-1}\big(\mathcal{W}(p)^{A_m}\big):=A_{0}\big(\mathcal{W}(p)^{A_m}\big).   F^{(m)}.
\end{gather*}

As in Lemma~\ref{uu} we have:
\begin{Lemma}
In $A(\mathcal{W}(p)^{A_m})$, we have $[E^{(m)}]^2=[F^{(m)}]^2=0$, $[E^{(m)}]*[F^{(m)}]\in
A_{0}(\mathcal{W}(p)^{A_m})$, and $[E^{(m)}]*[F^{(m)}]+[F^{(m)}]*[E^{(m)}]=g_p([\omega])$ for some fixed
$g_p(x)\in{\mathbb C}[x]$.

In particular,
\begin{gather*}
A\big(\mathcal{W}(p)^{A_m}\big):=A_{-1}\big(\mathcal{W}(p)^{A_m}\big) \oplus A_{0}\big(\mathcal{W}(p)^{A_m}\big) \oplus
A_{1}\big(\mathcal{W}(p)^{A_m}\big).
\end{gather*}
\end{Lemma}

\begin{Lemma}
\begin{gather*}
[H]*\big[F^{(m)}\big]=u_p([\omega]) * \big[F^{(m)}\big],
\end{gather*}
where the polynomial $u_p$ with $\deg u_p \leq p-1$ is uniquely determined by $u_p(h_{i,m+1})={- mp-1 + i \choose 2p-1}$, for $1\leq i\leq p$.
Similarly, we have
\begin{gather*}
[H]*\big[E^{(m)}\big]=-u_p([\omega]) * \big[E^{(m)}\big].
\end{gather*}
\end{Lemma}
\begin{proof}
The f\/irst assertion follows directly by evaluating this relation on lowest weight spaces of irreducible
$\mathcal{W}(p)^{A_m}$-modules.
The second assertion follows by applying the automorphism $\Psi \in \Aut(\mathcal{W}(p)^{A_m})$ such that
\begin{gather*}
\Psi \big(F^{(m)}\big)=E^{(m)}, \qquad  \Psi(H)=- H. \tag*{\qed}
\end{gather*}
\renewcommand{\qed}{}
\end{proof}

Similarly, we have
\begin{gather*}
\big[H \circ F^{(m)}\big]=v_p([\omega]) * \big[F^{(m)}\big]
\end{gather*}
with $\deg v_p \leq p$.
By evaluating this relation on lowest components of modules $\Lambda(i)^{-}$, we get that
\begin{gather*}
v_p(x)=k_{p} \prod\limits_{i=1}^{p} (x-h_{i,m+1}).
\end{gather*}
It follows from Conjecture~\ref{conjecutreA2} that $k_{p}$ is nonzero.
Hence we have
\begin{Proposition}\label{hefg}
Assume that Conjecture~{\rm \ref{conjecutreA2}} holds.
Then in $A(\mathcal{W}(p)^{A_m})$, we have
\begin{gather*}
\prod\limits_{i=1}^{p} ([\omega]-h_{i,m+1})\big[F^{(m)}\big]=0
\qquad \text{and} \qquad
\prod\limits_{i=1}^{p} ([\omega]-h_{i,m+1})\big[E^{(m)}\big]=0.
\end{gather*}
In particular,
\begin{gather*}
\dim A_{i} \big(\mathcal{W}(p)^{A_m}\big) \le p
\qquad
\text{for}
\quad
i=-1,1.
\end{gather*}
\end{Proposition}

Now we shall calculate the upper bound for $\dim A_0 (\mathcal{W}(p)^{A_m})$.

By the proof of~\cite[Lemma 2.1.3]{z},
\begin{gather*}
E^{(m)} \tilde{\circ} F^{(m)} -F^{(m)} \tilde{\circ} E^{(m)}=\Res_x \frac{(1+x)^{m^2p+mp-m}} {x^2}
Y(e^{-m\alpha},x) Q^{2m}e^{-m\alpha},
\\
E^{(m)} \tilde{\circ}_{4} F^{(m)} -F^{(m)} \tilde{\circ}_{4} E^{(m)}=\Res_x \frac{(1+x)^{m^2p+mp-m} (2+2x+x^{2})}
{x^4} Y(e^{-m\alpha},x) Q^{2m}e^{-m\alpha}.
\end{gather*}

On the other hand, by evaluating this relation on lowest components of $\mathcal{W}(p)^{A_m}$-modules we get

\begin{Lemma}
\begin{gather*}
  {\rm Res}_{x_{0},x_{1},\dots, x_{2m}}\frac{(1+x_{0})^{m^2p+mp-m(t+1)}} {x_{0}^{-2m^{2}p+2}(x_{1}\cdots x_{2m})^{2mp}}
\prod\limits_{i=1}^{2m} (1+ x_{i})^{t} \prod\limits_{1\leq i<j\leq
2m}(x_{i}-x_{j})^{2p}\prod\limits_{i=1}^{2m}(x_{i}-x_{0})^{-2mp}
\\
\qquad
=f_{p,m}(t){t-(m+1)p \choose p -1}{t+mp \choose p -1} \prod\limits_{i=-m^2p}^{m^2p} (t+1-p-\tfrac{i}{m}),
\\
  {\rm Res}_{x_{0},x_{1},\dots, x_{2m}}\frac{(1+x_{0})^{m^2p+mp-m}(2+2x_{0}+x_{0}^{2})} {x_{0}^{-2m^{2}p+4}(x_{1}\cdots
x_{2m})^{2mp}} (1+ x_{0})^{-mt}(1+ x_{1})^{t}\cdots(1+ x_{2m})^{t}
\\
\qquad
\phantom{=}
\times\prod\limits_{1\leq i<j\leq 2m}(x_{i}-x_{j})^{2p}\prod\limits_{i=1}^{2m}(x_{i}-x_{0})^{-2mp}
\\
\qquad
=\tilde{f}_{p,m}(t){t-(m+1)p \choose p -1}{t+mp \choose p -1} \prod\limits_{i=-m^2p}^{m^2p} (t+1-p-\tfrac{i}{m}),
\end{gather*}
where $f_{p,m}(t), \tilde{f}_p(t)\in {\mathbb C}[t]$, and
\begin{gather*}
\deg {f}_{p,m} \le 2(m-1)(p-1),
\qquad
\deg \tilde{f}_{p,m} \le 2(m-1)(p-1)+2,
\end{gather*}
\end{Lemma}

and use it to obtain (after we switch to $[\omega]$ polynomials)
\begin{gather*}
\big[E^{(m)} \tilde{\circ} F^{(m)} -F^{(m)} \tilde{\circ} E^{(m)}\big]=\big[F^{(m)} \circ E^{(m)}\big]
\\
\qquad
=s_{p,m}([\omega]) \prod\limits_{i=2p}^{mp+2p-1} ([\omega]- h_{i,1}) \prod\limits_{1\leq i\leq m^2p, m\nmid i}
\bigg([\omega]-h_{p+\tfrac{i}{m},1}\bigg)* [H],
\end{gather*}
and
\begin{gather*}
\big[E^{(m)} \tilde{\circ}_{4} F^{(m)} -F^{(m)} \tilde{\circ}_{4} E^{(m)}\big]
\\
\qquad{}
=\tilde{s}_{p,m}([\omega]) \prod\limits_{i=2p}^{mp+2p-1} ([\omega]- h_{i,1}) \prod\limits_{1\leq i\leq m^2p, m\nmid i}
\bigg([\omega]-h_{p+\tfrac{i}{m},1}\bigg)* [H],
\end{gather*}
where $s_{p,m}(x), \tilde{s}_{p,m}(x)\in {\mathbb C}[x]$, $\deg s_{p,m} \le (m-1)(p-1)$, $\deg \tilde{s}_{p,m} \le (m-1)(p-1)+1$.
Observe that~$s$-polynomials are of half degree of~$f$-polynomials.

\begin{Conjecture}
\label{conjectureA3}
Polynomials $f_{p,m}(t)$ and $\tilde{f}_{p,m}$ are relatively prime.
\end{Conjecture}

The next table gives evidence for the conjecture
\begin{center}\small
\begin{tabular}{|@{\,\,}c@{\,\,}|@{\,\,}c@{\,\,}|@{\,\,}c@{\,\,}|}
\hline
$(m,p)$ & $f_{p,m}(t)$ & $ \tilde{f}_{p,m}(t) $\tsep{2pt}\bsep{1pt}
\\
\hline
$(2,1) $& $ 1 $ & $4 t^2+85 $\tsep{1pt}
\\
\hline
$(2,2) $& $ 17 t^2-34 t+224 $ & $4 t^4-16 t^3+121 t^2-210 t+1568 $\tsep{1pt}
\\
\hline
$(2,3) $ & $ 233 t^4-1864 t^3+6539 t^2-11244 t+216216 $ & $764 t^6-9168 t^5+69807 t^4-313976 t^3$\tsep{1pt}
\\
& & ${}+895137 t^2-1459908 t+34378344 $
\\
\hline
$(2,4) $ & $ 811 t^6-14598 t^5+128467 t^4-665724 t^3 $ & $ 116908 t^8-2805792 t^7+32179967 t^6 $\tsep{1pt}
\\
$ $ & $ {} +1401172 t^2+422832 t+168030720 $ & $ -225709614 t^5+1213477115 t^4-5261506044 t^3 $
\\
& & $ {}+9465801460 t^2+10976862000 t+1491272640000 $
\\
\hline
$(3,1) $& $ 1 $ & $9 t^2+320 $\tsep{1pt}
\\
\hline
$(3,2) $& $ 26141 t^4-104564 t^3-576380 t^2$ & $4977 t^6-29862 t^5+495920 t^4-1784600 t^3$\tsep{1pt}
\\
& ${}+1361888 t-3720960$ & ${}-13382432 t^2+30254432 t-84341760 $
\\
\hline
\end{tabular}
\end{center}
If we assume that Conjecture~\ref{conjectureA3} holds, then $s_{p,m}(x)$ and $\tilde{s}_{p,m}(x)$ are relatively
prime.

\begin{Proposition}
\label{efg}
Assume that Conjecture~{\rm \ref{conjectureA3}} holds.
Then in $A(\mathcal{W}(p)^{A_m})$, we have
\begin{gather*}
(1)\quad \prod\limits_{i=2p}^{mp+2p-1} ([\omega]- h_{i,1}) \prod\limits_{1\leq i\leq m^2p, m\nmid i}
\bigg([\omega]-h_{p+\tfrac{i}{m},1}\bigg)* [H]=0,
\\
(2)\quad \prod\limits_{i=1}^{p}([\omega]- h_{i,1})\prod\limits_{i=1}^{p-1} ([\omega]- h_{mp+p+i,1}) \prod\limits_{i=1}^{m^2p}
\bigg([\omega]-h_{p+\tfrac{i}{m},1}\bigg)=0.
\end{gather*}
In particular,
\begin{gather*}
\dim A_{\omega} \big(\mathcal{W}(p)^{A_m}\big) \le \big(m^2 + 2\big)p -1,
\qquad
\dim A_{0} \big(\mathcal{W}(p)^{A_m}\big) \le \big(2 m^2 p + 2\big) p -1.
\end{gather*}
\end{Proposition}

\begin{proof}
The f\/irst assertion follows from the arguments which we explained above.
The second assertion follows by multiplying f\/irst identity by $[H]$.
\end{proof}

\begin{Theorem}
Assume that Conjectures~{\rm \ref{conjectureA3}} and~{\rm \ref{conjecutreA2}} $($or alternatively Conjectures~{\rm \ref{conjectureA3}}
and~{\rm \ref{conjectureA1})} hold.
Then
\begin{enumerate}\itemsep=0pt
\item[$(1)$] the above two tables give a~complete list of irreducible $\mathcal{W}(p)^{A_m}$-modules,
in particular, $\mathcal{W}(p)^{A_m}$ has $2m^2 p$ non-isomorphic irreducible modules,
\item[$(2)$] Zhu's algebra $A(\mathcal{W}(p)^{A_m})$ is of dimension $(2m^2 + 4)p-1$,
\item[$(3)$] the center of Zhu's algebra $A(\mathcal{W}(p)^{A_m})$ is $A_{\omega}(\mathcal{W}(p)^{A_m})$ and it has
dimension $(m^2 +2) p -1$.
\end{enumerate}
\end{Theorem}

\begin{proof}
It suf\/f\/ices to prove the second assertion.
By using Propositions~\ref{hefg} and~\ref{efg} we get
\begin{gather*}
\dim A\big(\mathcal{W}(p)^{A_m}\big)=\dim A_{-1} \big(\mathcal{W}(p)^{A_m}\big) + \dim A_0\big(\mathcal{W}(p)^{A_m}\big)
+ \dim A_1\big(\mathcal{W}(p)^{A_m}\big)
\\
\phantom{\dim A\big(\mathcal{W}(p)^{A_m}\big)}{}
 \le \big(2 m^2 p + 2\big) p -1 + 2p=\big(2 m^2 + 4\big) p-1.
\end{gather*}
Now assertion follows from Lemma~\ref{kriterij}.
\end{proof}

\begin{Remark}
In fact the classif\/ication of irreducible $\mathcal{W}(p)^{A_m}$-modules can also be derived from the classif\/ication of
irreducible $\mathcal{W}(p)^{D_m}$-modules by the general properties of orbifold vertex operator algebra.
\end{Remark}

\appendix

\section{Several conjectures about constant term identities}\label{ct-section}

In this section we list several constant term identities, although strictly speaking they are written as residues.
Let's start by recalling our identity from~\cite{ALM}:

\begin{Conjecture}[constant term identity of type~$A_m$,~I]
\label{conjectureA1}
\begin{gather*}
{\rm Res}_{x_0,x_1,\dots,x_{m+1}}\frac{(1+x_0)^{2p-1-t} \prod\limits_{i=1}^{m+1} (1+x_i)^t}{x_0^{2+2p} (x_1 \cdots
x_{m+1})^{2mp}}  \prod\limits_{i=1}^{m+1} \left(1-\frac{x_i}{x_0}\right)^{-2p} \prod\limits_{1 \leq i < j \leq m+1}(x_i-x_j)^{2p}
\\
\qquad
=\lambda_{p,m} {t+ mp \choose 2(m+1)p-1} \prod\limits_{i=1}^{m-1} {t+(i-1)p \choose 2ip-1},
\end{gather*}
where $\lambda_{p,m} \neq 0$.
\end{Conjecture}

In Sections~\ref{Section2} and~\ref{class-A}, we need the following two conjectures:

\begin{Conjecture}[constant term identity of type $A_m$,~II]\label{conjecutreA2}
\begin{gather*}
{\rm Res}_{x_{0},x_{1},\dots, x_{m+1}}\frac{(1+x_{0})^{m^2p+mp-m(t+1)}} {x_{0}^{-2mp+2}(x_{1}\cdots
x_{m+1})^{2mp}}\prod\limits_{i=1}^{m+1}(1+ x_{i})^{t} (x_{0}-x_{m+1})^{-2mp}
\\
\qquad
\phantom{=}{}
\times\prod\limits_{1\leq i<j\leq m+1}(x_{i}-x_{j})^{2p}\prod\limits_{i=1}^{m}(x_{i}-x_{0})^{-2mp}
\\
\qquad
=-2 {\rm Res}_{x_{0},x_{1},\dots, x_{m+1}}\frac{(1+x_{0})^{m^2p+mp-m(t+1)}} {x_{0}^{-2mp+3}(x_{1}\cdots
x_{m+1})^{2mp}}\prod\limits_{i=1}^{m+1}(1+ x_{i})^{t}   (x_{0}-x_{m+1})^{-2mp}
\\
\qquad
\phantom{=}
\times\prod\limits_{1\leq i<j\leq m+1}(x_{i}-x_{j})^{2p}\prod\limits_{i=1}^{m}(x_{i}-x_{0})^{-2mp}
\\
\qquad{}
=\mu_{p,m} {t-(m+1)p+1 \choose p} {t+mp \choose p} \prod\limits_{l=0}^{m-1}{t+pl \choose (m+1)p-1},
\end{gather*}
where $\mu_{p,m}$ is a~nonzero constant.
\end{Conjecture}
We verif\/ied this conjecture using Mathematica 9.0 for $m=2$, $p\leq12$; $m=3$, $p\leq6$; $m=4$, $p\leq4$.

\begin{Conjecture}[constant term identity of type $D_m$, $m>2$]
\label{conjecutreDm}
\begin{gather*}
{\rm Res}_{x_{0},x_{1},\dots, x_{m+2}}\frac{(1+x_{0})^{m^2p+mp-(m+1)t}} {x_{0}^{-4mp+3}(x_{1}\cdots
x_{m+2})^{4p}}\prod\limits_{i=1}^{m+2}(1+ x_{i})^{t}
\\
\qquad
\phantom{=}{}
\times(x_{0}-x_{m+1})^{-2mp}(x_{0}-x_{m+2})^{-2mp}\prod\limits_{1\leq i<j\leq
m+2}(x_{i}-x_{j})^{2p}\prod\limits_{i=1}^{m}(x_{i}-x_{0})^{-2mp}
\\
\qquad
=-{\rm Res}_{x_{0},x_{1},\dots, x_{m+2}}\frac{(1+x_{0})^{m^2p+mp-(m+1)t}} {x_{0}^{-4mp+4}(x_{1}\cdots
x_{m+2})^{4p}}\prod\limits_{i=1}^{m+2}(1+ x_{i})^{t}
\\
\qquad
\phantom{=}{}
\times(x_{0}-x_{m+1})^{-2mp}(x_{0}-x_{m+2})^{-2mp}\prod\limits_{1\leq i<j\leq
m+2}(x_{i}-x_{j})^{2p}\prod\limits_{i=1}^{m}(x_{i}-x_{0})^{-2mp}
\\
\qquad
=\alpha_{p,m} {t+p+\tfrac{1}{2} \choose 2p} {t-p+ \tfrac{1}{2} \choose 2p} {t-(m+1)p+1 \choose p} {t+mp \choose p}
\prod\limits_{l=0}^{m-1}{t+pl \choose (m+1)p-1},
\end{gather*}
where $m\geq3$ and $\alpha_{p,m}$ is a~nonzero constant.

\end{Conjecture}
We verif\/ied this conjecture using Mathematica~9.0 for $m=3$, $p\leq4$; $m=4$, $p\leq2$.

Note that the Conjecture~\ref{conjecutreDm} does not hold for $m=2$.
Instead we have the following conjecture (cf.~\cite{ALM1}).

\begin{Conjecture}[\protect{\cite[Conjecture~7.6]{ALM1}}, constant term identity of type~$D_2$, I)]
\label{conjecutreD2}
\begin{gather*}
{\rm Res}_{x_{0},x_{1},x_{2},x_{3},x_{4}}\frac{(1+x_{0})^{6p-2-2t}} {x_{0}^{-8p+3}(x_{1}x_{2}x_{3} x_{4})^{4p}}(1+
x_{1})^{t}(1+ x_{2})^{t}(1+ x_{3})^{t}(1+ x_{4})^{t}
\\
\qquad
\phantom{=}
\times(x_{0}-x_{1})^{-4p}(x_{0}-x_{2})^{-4p}(x_{3}-x_{0})^{-4p}\prod\limits_{1\leq i<j\leq
4}(x_{i}-x_{j})^{2p}\frac{\partial_{x_{0}}^{4p-1}}{(4p-1)!} x_{4}^{-1} \delta\left(\frac{x_{0}}{x_{4}}\right)
\\
=A_p {t + p + 1/2 \choose 4p} {t + 2p \choose 4p -1} {t \choose 4p -1},
\end{gather*}
where $A_p$ is a~nonzero constant.
\end{Conjecture}

It is very interesting to note that Conjectures~\ref{conjecutreA2} and~\ref{conjecutreDm} have no common natural
gene\-ra\-li\-zation.
We close this paper with another beautiful constant term identity of type $D_2$ which we have verif\/ied for $p\leq6$:

\begin{Conjecture}[\textbf{Constant term identity of type~$D_2$,~II}]
\begin{gather*}
  {\rm Res}_{x_{0},x_{1},x_{2},x_{3},x_{4}}\frac{(1+x_{0})^{6p-2-2t}} {x_{0}^{-8p+2}(x_{1}x_{2}x_{3} x_{4})^{4p}}(1+
x_{1})^{t}(1+ x_{2})^{t}(1+ x_{3})^{t}(1+ x_{4})^{t}
\\
\qquad\phantom{=}{}\times  (x_{0}-x_{1})^{-4p}(x_{0}-x_{2})^{-4p}(x_{3}-x_{0})^{-4p}(x_{4}-x_{0})^{-4p}\prod\limits_{1\leq i<j\leq
4}(x_{i}-x_{j})^{2p}
\\
\qquad
=D_p {t+2p \choose 6p-1} {t+p \choose 4p-1} {t \choose 2p-1},
\end{gather*}
where $D_p$ is a~nonzero constant.
\end{Conjecture}

\begin{Remark}
In all conjectural identities we expect constants $\lambda_{p,m}$, $\mu_{p,m}$, $\alpha_{p,m}$, $A_p$ and $D_p$ to be
expressible in terms of quotients of binomial coef\/f\/icients which depend on~$m$ and~$p$ linearly.
This is further supported by our numerical calculations.
\end{Remark}

\subsection*{Acknowledgements}

We would like to thank the referees for their valuable comments.
D.A.\
is partially supported by the Croatian Science Foundation under the project~2634.
X.L.\
is partially supported by National Natural Science Foundation for young (no.~11401098).
A.M.\
is partially supported by a~Simons Foundation grant.

\pdfbookmark[1]{References}{ref}
\LastPageEnding


\begin{thebibliography}{99}
\footnotesize \itemsep=0pt

\bibitem{A-2003}
Adamovi{\'c} D., Classif\/ication of irreducible modules of certain subalgebras
  of free boson vertex algebra, \href{http://dx.doi.org/10.1016/j.jalgebra.2003.07.011}{\textit{J.~Algebra}} \textbf{270} (2003),
  115--132, \href{http://arxiv.org/abs/math.QA/0207155}{math.QA/0207155}.

\bibitem{ALM}
Adamovi{\'c} D., Lin X., Milas A., A{DE} subalgebras of the triplet vertex
  algebra {${\mathcal W}(p)$}: {$A$}-series, \href{http://dx.doi.org/10.1142/S0219199713500284}{\textit{Commun. Contemp. Math.}}
  \textbf{15} (2013), 1350028, 30~pages, \href{http://arxiv.org/abs/1212.5453}{arXiv:1212.5453}.

\bibitem{ALM1}
Adamovi{\'c} D., Lin X., Milas A., A{DE} subalgebras of the triplet vertex
  algebra {${\mathcal W}(p)$}: {$D$}-series, \href{http://dx.doi.org/10.1142/S0129167X14500013}{\textit{Internat.~J. Math.}}
  \textbf{25} (2014), 1450001, 34~pages, \href{http://arxiv.org/abs/1304.5711}{arXiv:1304.5711}.

\bibitem{am-1}
Adamovi{\'c} D., Milas A., Logarithmic intertwining operators and {${\mathcal
  W}(2,2p-1)$} algebras, \href{http://dx.doi.org/10.1063/1.2747725}{\textit{J.~Math. Phys.}} \textbf{48} (2007), 073503,
  20~pages, \href{http://arxiv.org/abs/math.QA/0702081}{math.QA/0702081}.

\bibitem{am1}
Adamovi{\'c} D., Milas A., On the triplet vertex algebra {${\mathcal W}(p)$},
  \href{http://dx.doi.org/10.1016/j.aim.2007.11.012}{\textit{Adv. Math.}} \textbf{217} (2008), 2664--2699, \href{http://arxiv.org/abs/0707.1857}{arXiv:0707.1857}.

\bibitem{am-sigma}
Adamovi{\'c} D., Milas A., The {$N=1$} triplet vertex operator superalgebras:
  twisted sector, \href{http://dx.doi.org/10.3842/SIGMA.2008.087}{\textit{SIGMA}} \textbf{4} (2008), 087, 24~pages,
  \href{http://arxiv.org/abs/0806.3560}{arXiv:0806.3560}.


\bibitem{am-imrn}
Adamovi{\'c} D., Milas A., On {$W$}-algebras associated to {$(2,p)$} minimal
  models and their representations, \href{http://dx.doi.org/10.1093/imrn/rnq016}{\textit{Int. Math. Res. Not.}} \textbf{2010}
  (2010), 3896--3934, \href{http://arxiv.org/abs/0908.4053}{arXiv:0908.4053}.

\bibitem{am-ja}
Adamovi{\'c} D., Milas A., On {${\mathcal W}$}-algebra extensions of {$(2,p)$}
  minimal models: {$p>3$}, \href{http://dx.doi.org/10.1016/j.jalgebra.2011.07.006}{\textit{J.~Algebra}} \textbf{344} (2011), 313--332,
  \href{http://arxiv.org/abs/1101.0803}{arXiv:1101.0803}.

\bibitem{am2}
Adamovi{\'c} D., Milas A., The structure of {Z}hu's algebras for certain
  {${\mathcal W}$}-algebras, \href{http://dx.doi.org/10.1016/j.aim.2011.05.007}{\textit{Adv. Math.}} \textbf{227} (2011),
  2425--2456, \href{http://arxiv.org/abs/1006.5134}{arXiv:1006.5134}.

\bibitem{dg}
Dong C., Griess Jr. R.L., Rank one lattice type vertex operator algebras and
  their automorphism groups, \href{http://dx.doi.org/10.1006/jabr.1998.7498}{\textit{J.~Algebra}} \textbf{208} (1998), 262--275,
  \href{http://arxiv.org/abs/q-alg/9710017}{q-alg/9710017}.

\bibitem{fw}
Forrester P.J., Warnaar S.O., The importance of the {S}elberg integral,
  \href{http://dx.doi.org/10.1090/S0273-0979-08-01221-4}{\textit{Bull. Amer. Math. Soc.}} \textbf{45} (2008), 489--534,
  \href{http://arxiv.org/abs/0710.3981}{arXiv:0710.3981}.

\bibitem{flm}
Frenkel I., Lepowsky J., Meurman A., Vertex operator algebras and the
  {M}onster, \textit{Pure and Applied Mathematics}, Vol.~134, Academic Press,
  Inc., Boston, MA, 1988.

\bibitem{gin}
Ginsparg P., Curiosities at {$c=1$}, \href{http://dx.doi.org/10.1016/0550-3213(88)90249-0}{\textit{Nuclear Phys.~B}} \textbf{295}
  (1988), 153--170.

\bibitem{ll}
Lepowsky J., Li H., Introduction to vertex operator algebras and their
  representations, \href{http://dx.doi.org/10.1007/978-0-8176-8186-9}{\textit{Progress in Mathematics}}, Vol.~227, Birkh\"auser
  Boston, Inc., Boston, MA, 2004.

\bibitem{z}
Zhu Y., Modular invariance of characters of vertex operator algebras,
  \href{http://dx.doi.org/10.1090/S0894-0347-96-00182-8}{\textit{J.~Amer. Math. Soc.}} \textbf{9} (1996), 237--302.

\end{thebibliography}
\end{document}